\documentclass[preprint,authoryear,12pt]{elsarticle}

\usepackage{amssymb}

\usepackage{amsmath,amsfonts}
\usepackage{times,mathptm}
\usepackage{cite}
\usepackage{amsmath,amsfonts,epsfig,subfigure}
\usepackage[latin1]{inputenc}

\newcommand{\sgn}{\text{sgn}}
\newcommand{\e}{\text{e}}
\DeclareMathAlphabet{\mathbi}{OT1}{cmr}{bx}{it}
\SetMathAlphabet{\mathbi}{normal}{OT1}{cmr}{bx}{it}

\newcommand{\C}{\mathbf{C}}

\newcommand{\tr}{\text{tr}\,}

\newtheorem{definition}{Definition}

\newtheorem{remark}[definition]{Remark}

\begin{document}

\begin{frontmatter}



\title{Complex dynamics in hard oscillators: the influence of
constant inputs}


\author[delen]{V. Lanza\corref{vale}}
\ead{valentina.lanza@polito.it}

\author[difis]{L. Ponta}
\ead{linda.ponta@polito.it}

\author[delen]{M. Bonnin}
\ead{michele.bonnin@polito.it}

\author[delen]{F. Corinto}
\ead{fernando.corinto@polito.it}

\address[delen]{Department of Electronics, Politecnico di Torino, Italy}
\address[difis]{Department of Physics, Politecnico di Torino, Italy}

\cortext[vale]{Corresponding author}

\begin{abstract}
Systems with the coexistence of different stable
attractors are widely exploited in systems biology in order to
suitably model the differentiating processes arising in living
cells. In order to describe genetic regulatory networks several
deterministic models based on systems of nonlinear ordinary
differential equations have been proposed.

Few studies have been developed to characterize how either an external
input or the coupling can drive systems with different coexisting
states. For the sake of simplicity, in this manuscript we focus on systems
belonging to the class of radial isochron
clocks that exhibits hard excitation, in order to investigate their complex
dynamics, local and global bifurcations arising in presence
of constant external inputs. In particular the occurrence of
saddle node on limit cycle bifurcations is detected.

%

\end{abstract}


\end{frontmatter}



\section{Introduction}


Being the most diffuse formalism to model dynamical systems in
science and engineering, ordinary differential equations (ODEs)
have been widely used in systems biology as well. In particular,
nonlinear dynamical systems that can exhibit coexisting stable
attractors are considered of universal importance, since they
allow physiologists to accurately model cell differentiation
processes. One of the most interesting cases, that goes by the name of of hard
excitation, is the concurrence between oscillatory and non-oscillatory states
\citep{liu,sontag2008,minorsky}.


It is well known the
external environment (such as light and temperature) or the substrate
synthesis/injection rate play an important role in the dynamics of
molecular reactions \citep{bastin}. In order to properly model these effects,
external forcing terms have to be taken into account in
the mathematical modeling \citep{leloup}.
Moreover, the processes inside the cell are localized
in different spatial domains, while exchanges of chemical species take place in
the common extracellular medium. Therefore, it is more
appropriate to consider reaction diffusion or multi-compartmental equations
\citep{sontag2008}. Actually, few studies \citep{goldbeter,goldbeter_book}
have
been developed to characterize how an external input or coupling effects can
drive systems with different coexisting states.



Typical systems that exhibit such a dynamical
behavior are represented by the so-called Cyclic
Negative Feedback Systems, that arise in a variety
of mathematical bio-inpired models, from cellular signal pathways
\citep{kholodenko,liu} to gene regulatory networks
\citep{othmer,dejong,elowitz}. Cyclic Negative Feedback systems (CNF systems)
are
described by
the following set of nonlinear differential equations
\citep{arcak2008}:
\begin{equation}
\label{eq:Nfcs-gen}
 \begin{cases}
  \dot x_1  = -g_1(x_1)+f_n(x_n)\nonumber\\
\dot x_k  = -g_k(x_k)+f_{k-1}(x_{k-1}), \qquad 2\leq k\leq n
 \end{cases}
\end{equation}
where $x_k$ only assumes positive values, $g_p(\cdot)$
($p=1,\dots, n$) and $f_q(\cdot)$ ($q=1,\dots, n-1$) are increasing
functions, while $f_n(\cdot)$ is a decreasing function. The
first-order time derivative of $x_k$ is denoted by $\dot x_k$.
The variables $x_k$ can represent the concentrations of certain
molecules (e.g. mRNA, proteins) in the cell. In addition, the function
$f_n(\cdot)$
suitably reproduces the inhibition of $x_1$ by the chemical
product $x_n$. It is easy to deduce from the previous expression that each variable
$x_k(t)$ is activated by its previous neighbor $x_{k-1}(t)$, except for $x_1(t)$
that is repressed by $x_n(t).$ It can be proved \citep{arcak2008} that CNF systems
exhibit a global asymptotic stable equilibrium $\mathbf
{x}^\ast=(x_1^\ast,\dots, x_n^\ast),$ provided the secant
criterion \citep{thron} is satisfied. 
Furthermore, for these type of systems the Poincaré-Bendixson
Theorem holds \citep{malletparet}, that is the coexistence between
stable equilibria and periodic orbits is allowed.


In \citep{lanza_ijcnn09}, we have analyzed the effect of coupling on
arrays of diffusively coupled third order CNF systems. We have shown that CNF
arrays with diffusive couplings that are constant and local (i.e. they involve
only the two nearest neighbors of each CNF system) are potentially equivalent to
nonlinear networks whose elements are fully connected (i.e. each subsystem is
linked to all the others).
Moreover, we have shown that already in a two-compartment version of CNF systems
new dynamics, such as global periodic oscillations with
space-variant amplitude (e.g., discrete breathers-like patterns), can arise due
to the couplings.

One of the main drawbacks of CNF systems is that they can be of high order and
quite difficult to handle. Even for the network of third order systems
analyzed in \citep{lanza_ijcnn09}, it is matematically complicated to detect
and characterize all the new periodic
solutions to whom the coupling or an external input give rise.
Therefore, in order to investigate the emergence of such new
dynamics, we consider a nonlinear dynamical system that is easier
to handle and that for certain values of its parameters has
qualitatively the same dynamical behavior of a CNF system. First
of all, we focus on the single system and carry out a complete
analysis of both the local and global bifurcations arising in
presence of constant external inputs. In particular, we show the
occurrence of saddle node on limit cycle bifurcations and that in
presence of such bifurcations our system behaves as a relaxation
oscillator.

This manuscript is structured as follows. In Section 2 we
introduce the model under study, i.e. a radial isochron clock with
hard excitation that is well known in literature as the normal
form for the Bautin bifurcation \citep{kuznetsov,izhikevich2001}.
Moreover, we investigate how this system changes its form due to
the presence of a constant external input. In Section 3 a complete
analysis of local and global bifurcations arising in presence of
constant external inputs is carried out. In particular, we show
that the periodic solutions of our system can disappear through
saddle node on limit cycle bifurcations. Section 4 is devoted to
conclusions.

\section{The case of radial isochron clocks}

Let us consider a simpler model with respect to CNF systems but with the same
dynamics:
\begin{equation}\label{eq:systcomplex}
 \dot z= (\sigma_0+j\Omega_0)z+\sigma_1 |z|z+\sigma_2 |z|^2 z+I_0,
\end{equation}
where $z=r\e^{j\phi}\in \C$ is a complex variable, $\Omega_0>0$ and $\sigma_0,
\sigma_1, \sigma_2$ are real parameters. The term $I_0$ represents the
action of a constant external input and from now on we will
assume $I_0>0$ without any loss of generality\footnote{It is easy to notice that
the system in invariant under the transformation $z\to -z, I_0\to -I_0$.}.

\subsection*{Absence of constant external input}
In the $(r,\phi)$ coordinates, system
\eqref{eq:systcomplex} with $I_0=0$ can be recast as:
\begin{equation}\label{eq:singlesyst}
\begin{cases}
\dot r  = r(\sigma_0+\sigma_1r+\sigma_2r^2)\\
\dot \phi  = \Omega_0.
\end{cases}
\end{equation}
It is worth noting that \eqref{eq:singlesyst} belongs to the category of the
radial isochron clocks proposed by Winfree \citep{winfree}, and it has been
widely exploited as a simplified model for the
Hodgkin-Huxley neuron \citep{curran} and for studying the phenomenon of
cardiac
fibrillation \citep{depaor}.

\begin{figure}[t!]
\centering
\includegraphics[bb=32 66 426 359,width=65mm,clip]{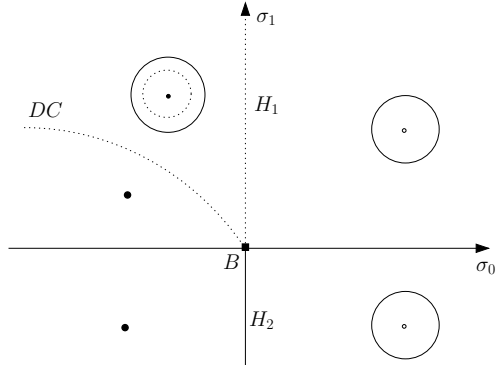}
\caption{Bifurcation diagram for the single radial isochron clock
\eqref{eq:singlesyst}. The different bifurcation curves are
labeled as follows: $H_1$ stands for subcritical Hopf bifurcation, $H_2$
for supercritical Hopf bifurcation, $DC$ for double limit cycle
bifurcation, and $B$ for Bautin bifurcation. The dotted line where a
$DC$ bifurcation occurs has equation $\sigma_1^2-\sigma_0
\sigma_2=0$.}\label{fig:bautin}
\end{figure}

This system can be easily investigated since the study of the periodic solutions
of \eqref{eq:singlesyst} reduces to the analysis of the equilibria of the first
equation of \eqref{eq:singlesyst}. It is easy to see that
\begin{itemize}
 \item a Hopf bifurcation takes place for
$\sigma_0=0$, which is subcritical for $\sigma_1>0$ and supercritical
otherwise;
\item a double limit cycle
bifurcation (also known as fold, or tangent, or saddle-node
bifurcation of limit cycles) occurs for $\sigma_1^2-4\sigma_0\sigma_2=0$,
and $\sigma_1\geq 0$. A double limit cycle bifurcation occurs when
a branch of stable periodic solutions and a branch of unstable
periodic solutions coalesce and obliterate each other at the
bifurcation point \citep{guckenheimer,kuznetsov};
\item this system can undergo a Bautin bifurcation, when a Hopf
bifurcation and a double limit cycle bifurcation occur simultaneously
\citep{kuznetsov,izhikevich2001}. In our case it happens for
$\sigma_0=\sigma_1=0$ but $\sigma_2\neq 0$. When $\sigma_2<0$, the Bautin
bifurcation is said to be supercritical, while for $\sigma_2>0$ it is
subcritical. In particular, for $\sigma_2<0$ the cycle with larger amplitude
is stable.
\end{itemize}
The complete bifurcation diagram is represented in Figure
\ref{fig:bautin}.

\begin{figure}[t!]
\centering
\includegraphics[bb=158 520 425 792,width=40mm,clip]{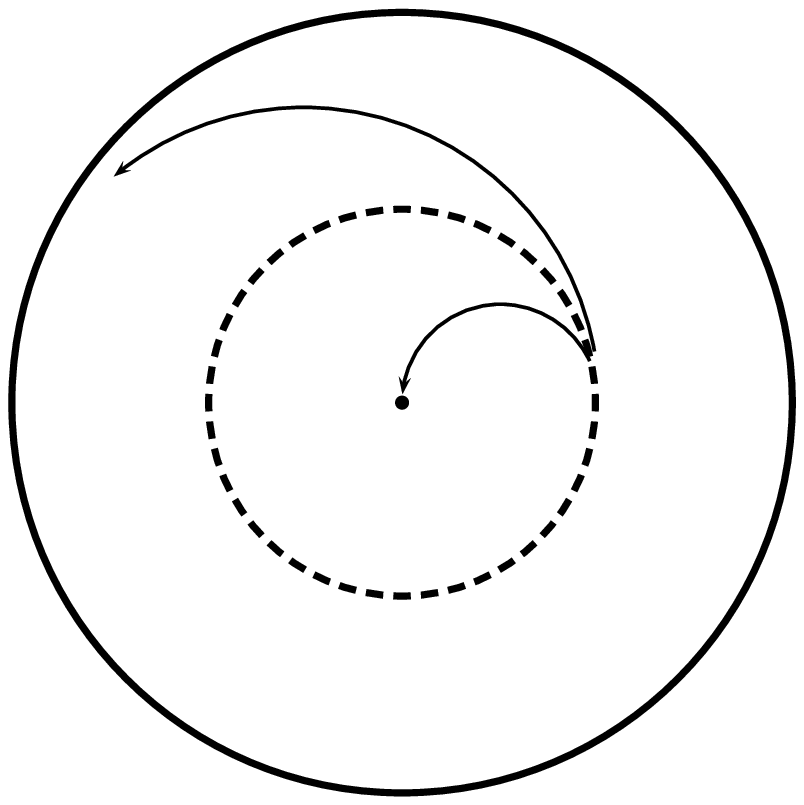}
\caption{Phase portrait of a system with hard excitation. The
stable equilibrium point and the stable limit cycle (solid line)
are separated by an unstable limit cycle (dashed line).}
\label{fig:phaseportrait}
\end{figure}

In the following, we are interested in studying a system that
presents a hard excitation behavior, thus we focus on a choice of
the parameters such that we have the coexistence of a stable
equilibrium point and a stable limit cycle:
\begin{equation}\label{eq:cond2cycles}
 \begin{cases}
  \sigma_1^2-4\sigma_0\sigma_2>0\\
  \sigma_1>0\\
  \sigma_2<0.
 \end{cases}
\end{equation}
In Figure \ref{fig:phaseportrait} the phase portrait of our system is
represented. We have a sort of concentric structure of cycles, where the origin
 is a stable equilibrium point, sorrounded by two closed curves alternatively
stable and unstable. Thus, it is easy to understand why a similar
configuration goes by the name of hard excitation
\citep{minorsky}: if the systems is in the steady state, then it
needs a strong perturbation to cross the separatrix and start to
oscillate.

\subsection*{Presence of external constant input}
In order to simplify the notation and reduce
the number
of control parameters, we rescale the variables involved in
\eqref{eq:systcomplex} in the following way:
$$(z,t) \to
\left(\sqrt{\frac{|\sigma_0|}{|\sigma_2|}}z,\frac{t}{|\sigma_0|}\right).$$
Thus, from \eqref{eq:systcomplex} we obtain the system
\begin{equation}\label{eq:systcomplex_input_risc}
 \dot z= (\sgn\, \sigma_0+j\Omega)z+\alpha |z|z+\sgn\, \sigma_2 |z|^2 z
+I,
\end{equation}
where $\sgn(\cdot)$ is the sign function and
\begin{equation}
 \Omega=\frac{\Omega_0}{|\sigma_0|}\qquad
\alpha=\frac{\sigma_1}{\sqrt{|\sigma_0 \sigma_2|}} \qquad
I=\sqrt{\left|
\frac{\sigma_2}{\sigma_0}\right|}\frac{I_0}{|\sigma_0|}.
\end{equation}
Conditions \eqref{eq:cond2cycles} imply $\sgn\, \sigma_0=-1,
\sgn\, \sigma_2=-1, \alpha>2$, and thus system
\eqref{eq:systcomplex_input_risc} becomes
\begin{equation}\label{eq:systcomplex_input_risc2}
 \dot z= (-1+j\Omega)z+\alpha |z|z- |z|^2 z
+I.
\end{equation}
Exploiting the cartesian coordinates ($z=x+jy$),
system \eqref{eq:systcomplex_input_risc2} can be
rewritten as
\begin{equation}\label{eq:systcart_input}
\begin{cases}
\dot x = (-1+\alpha\sqrt{x^2+y^2}-(x^2+y^2))x-\Omega y+I\\
\dot y = (-1+\alpha\sqrt{x^2+y^2}-(x^2+y^2))y+\Omega x,
\end{cases}
\end{equation}
or to simplify notation
\begin{equation}\label{eq:systcart_input2}
\begin{cases}
\dot x = g(x,y)x-\Omega y+I\\
\dot y = g(x,y)y+\Omega x,
\end{cases}
\end{equation}
 where $g(x,y)=(-1+\alpha\sqrt{x^2+y^2}-(x^2+y^2))$.

First of all, in order to characterize the dynamics of this
system, we are interested in finding all the equilibrium
configurations. Thus, we have to solve the following set of
nonlinear algebraic equations:
\begin{equation}\label{eq:systeq}
 \begin{cases}
  g(x,y)x-\Omega y+I  =0\\
g(x,y)y+\Omega x  =0.
 \end{cases}
\end{equation}
From the second equation we obtain
\begin{equation}\label{eq:g}
 g(x,y)=-\Omega \frac{x}{y},
\end{equation}
having assumed\footnote{The case $y=0$ can be treated separately and
is not interesting for our study. In fact, from \eqref{eq:systeq} we can
conclude that $y=0$ implies $x=0$ and $I=0$. Thus, a solution with $y=0$ can
be achieved only in absence of external input.} $y\neq 0$.
The substitution of expression \eqref{eq:g} in the first equation of
\eqref{eq:systeq} yields
\begin{equation}\label{eq:circum}
-\Omega\frac{x^2}{y}-\Omega y+I=0.
\end{equation}
The solutions of \eqref{eq:systcart_input2} are precisely the
intersection points between the circumference $\Gamma:
x^2+\left(y-\frac{I}{2\Omega}\right)^2=\left(\frac{I}{2\Omega}\right)^2$
and the curve $S:x=-\frac{1}{\Omega}y g(x,y)$. Moreover, from
\eqref{eq:circum} we can conclude that all the equilibrium
configurations have a positive $y-$component. Furthermore, since
from \eqref{eq:circum} we have $x^2+y^2=\frac{I}{\Omega}y,$ we can
notice that actually $S$ has the following expression:
$$S: x=-\frac{1}{\Omega}\left[-1+\alpha\sqrt{\frac{I}{\Omega}y}
-\frac{I}{\Omega}y\right ],$$ and we can derive that the
equilibrium points of system \eqref{eq:systcart_input2} are the
intersections between the curves:
\begin{equation}\label{eq:intersez}
 \begin{cases}
  \Gamma:  \,\,\,
x^2+\left(y-\frac{I}{2\Omega}\right)^2=\left(\frac{I}{2\Omega}\right)^2\\
S:  \,\,\,
x=-\frac{1}{\Omega}\left[-1+\alpha\sqrt{\frac{I}{\Omega}y}
-\frac{I}{\Omega}y\right ].
 \end{cases}
\end{equation}
Substituting the second equation of \eqref{eq:intersez} in the first one, and
introducing the new variable $\zeta=\sqrt{\frac{I}{\Omega}y}$, we finally
conclude that the equilibrium points of system \eqref{eq:systcart_input} are the
roots of the following polynomial:
\begin{equation}\label{eq:polyeq}
\zeta^6-2\alpha
\zeta^5+(\alpha^2+2)\zeta^4-2\alpha\zeta^3+(1+\Omega^2)\zeta^2-I^2=0.
\end{equation}
It is worth noting that, due to the definition of $\zeta$, we are
interested only in the real and positive roots of
\eqref{eq:polyeq}.

\begin{figure}[ht!]
 \begin{minipage}[l]{5.5cm}
   \centering
   \includegraphics[bb=81 242 507 592,width=55mm,clip]{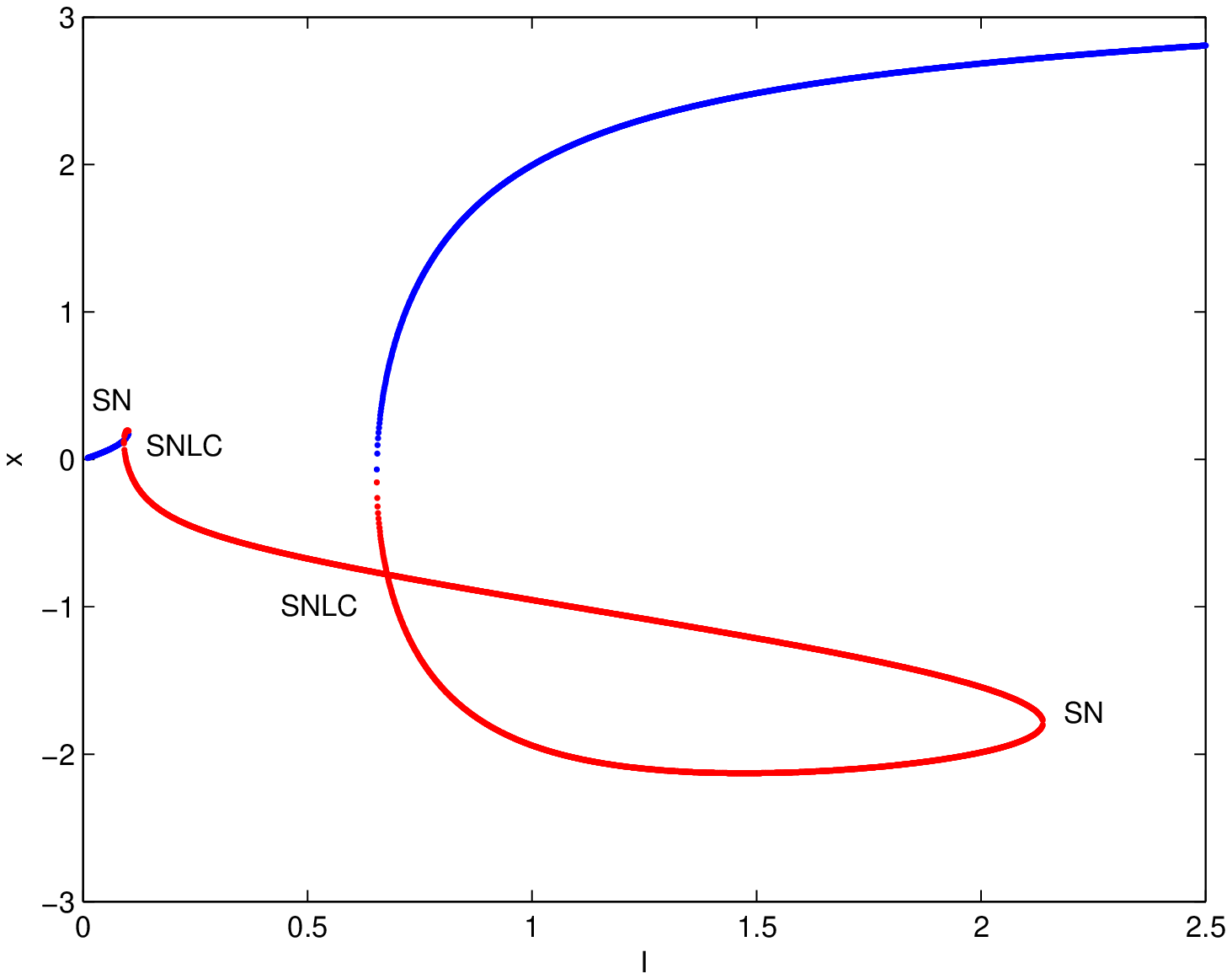}\\
   (a)
 \end{minipage}
\qquad
 \begin{minipage}[r]{5.5cm}
    \centering
    \includegraphics[bb=81 242 507 592,width=55mm,clip]{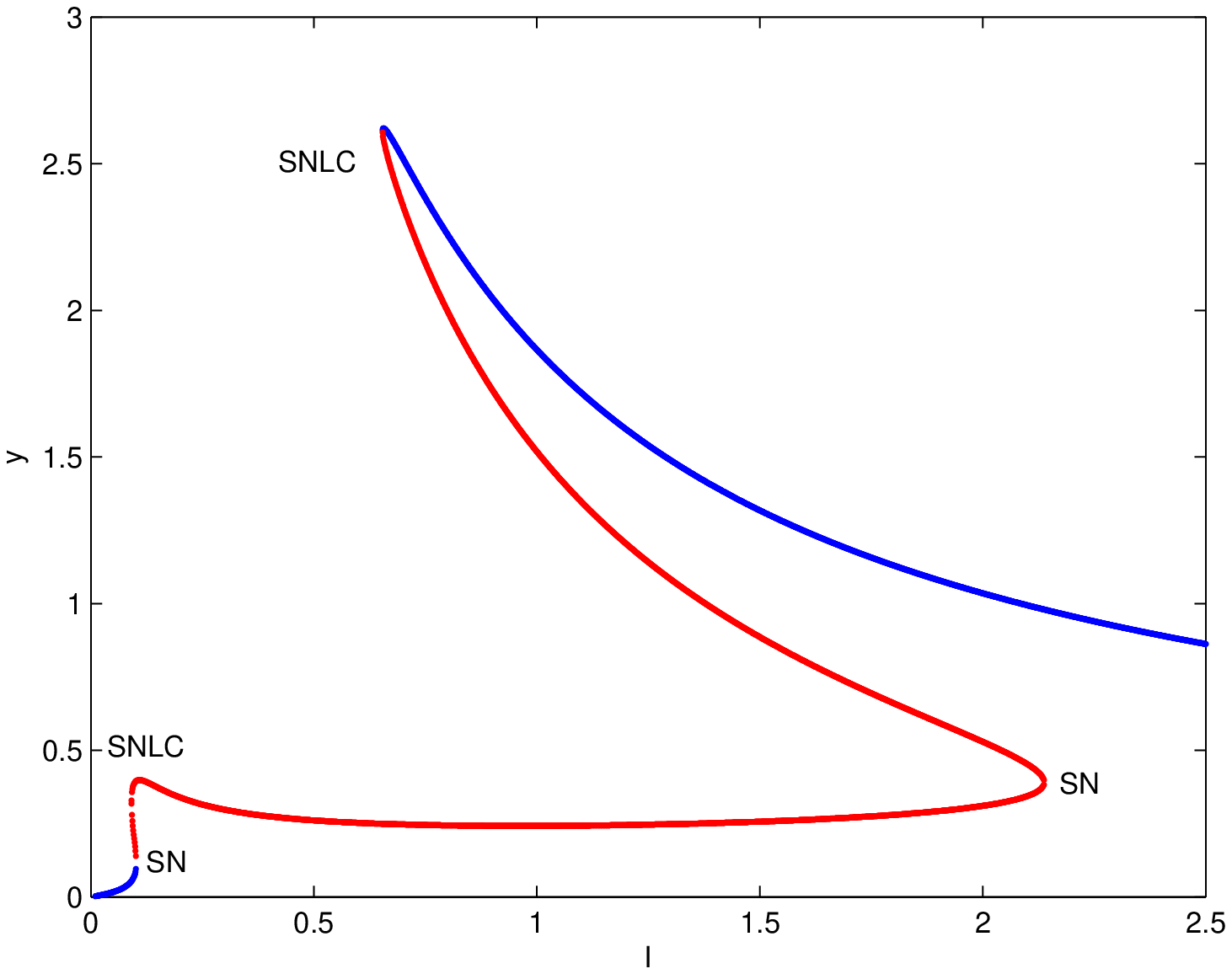}\\
    (b)
 \end{minipage}
 \\[0.3cm]
   \begin{minipage}[l]{5.5cm}
   \centering
   \includegraphics[bb=81 242 507 592,width=55mm,clip]{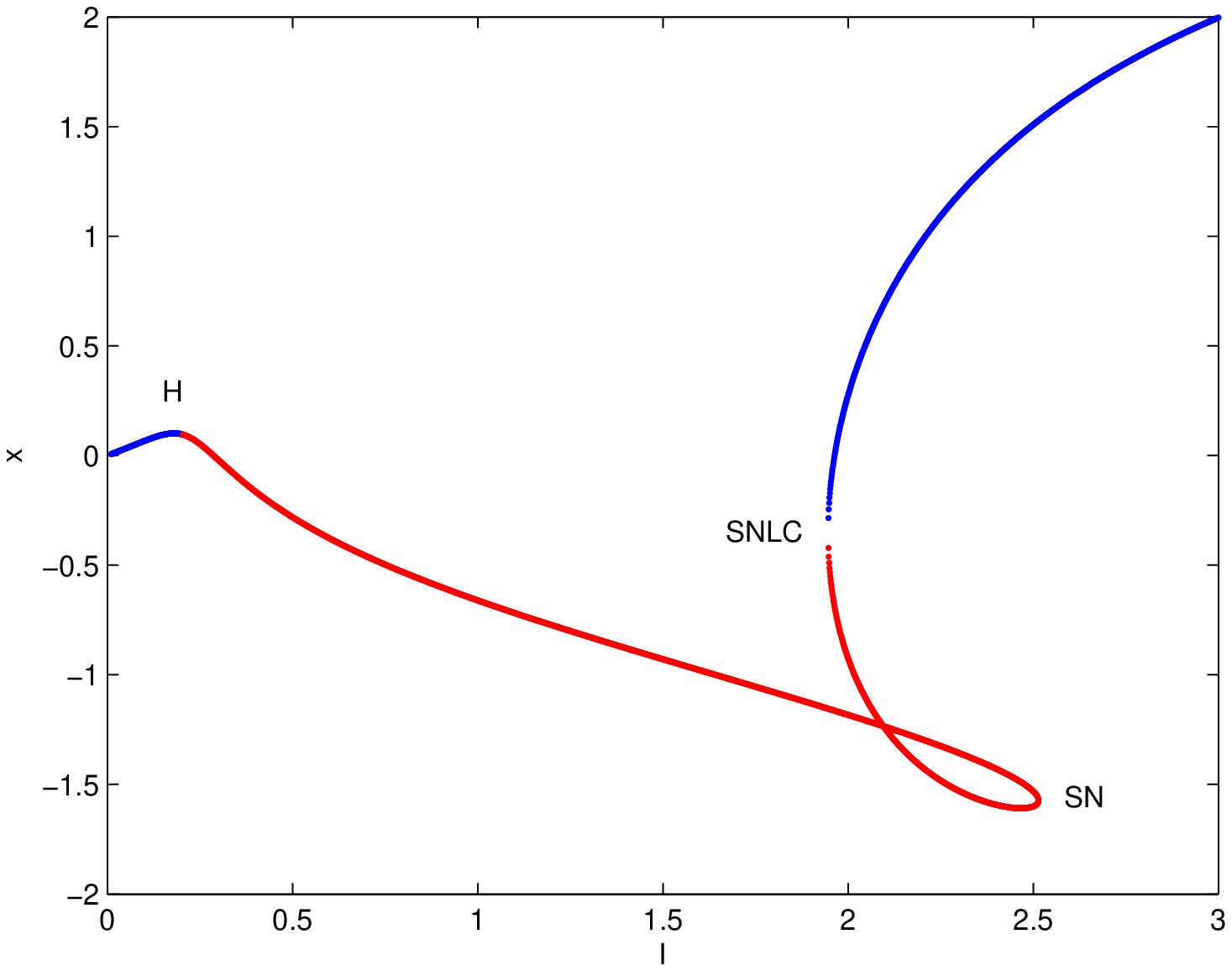}\\
   (c)
 \end{minipage}
\qquad
 \begin{minipage}[r]{5.5cm}
    \centering
    \includegraphics[bb=81 242 507 592,width=55mm,clip]{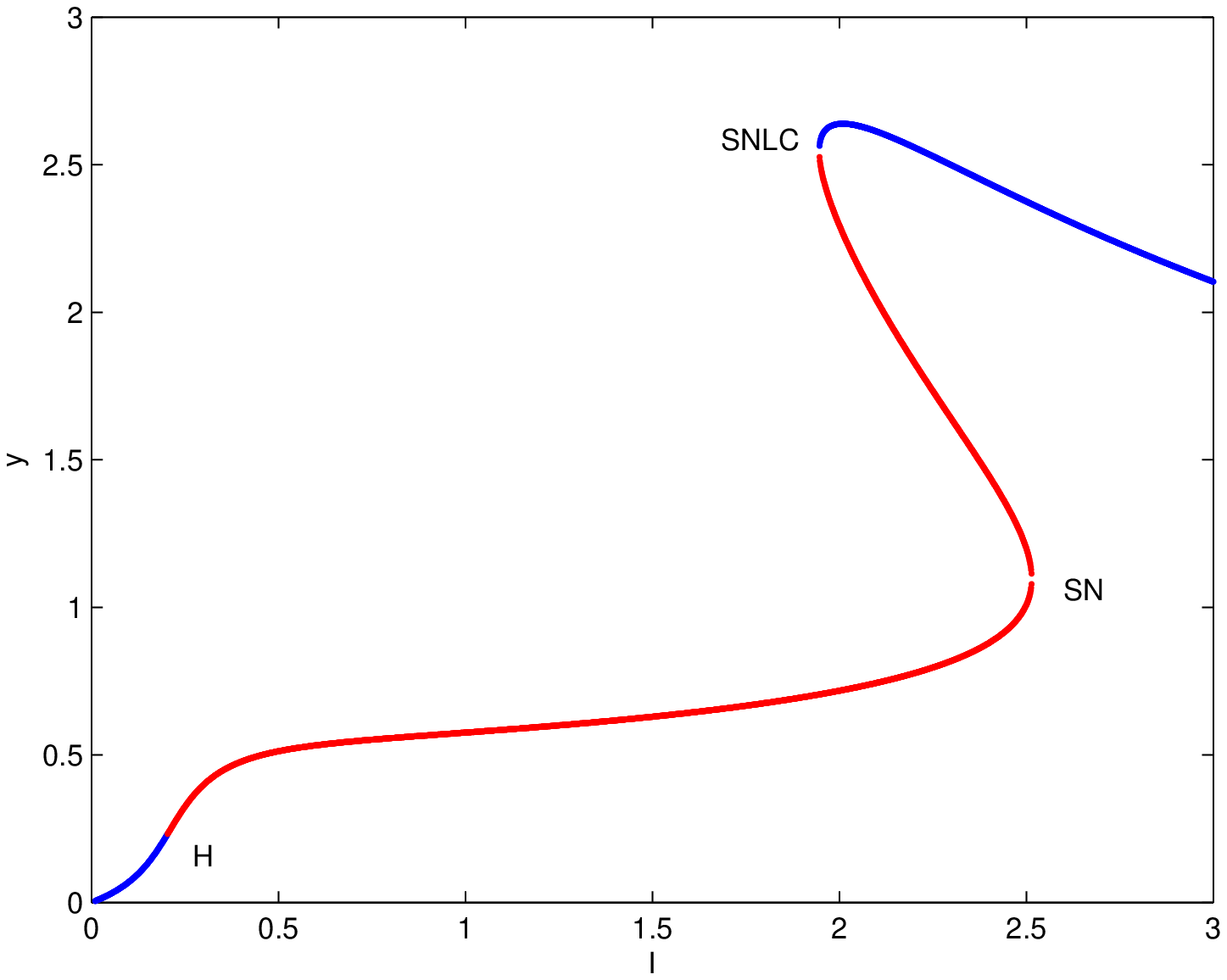}\\
    (d)
 \end{minipage}
 \\[0.3cm]
   \begin{minipage}[l]{5.5cm}
   \centering
   \includegraphics[bb=81 242 507 592,width=55mm,clip]{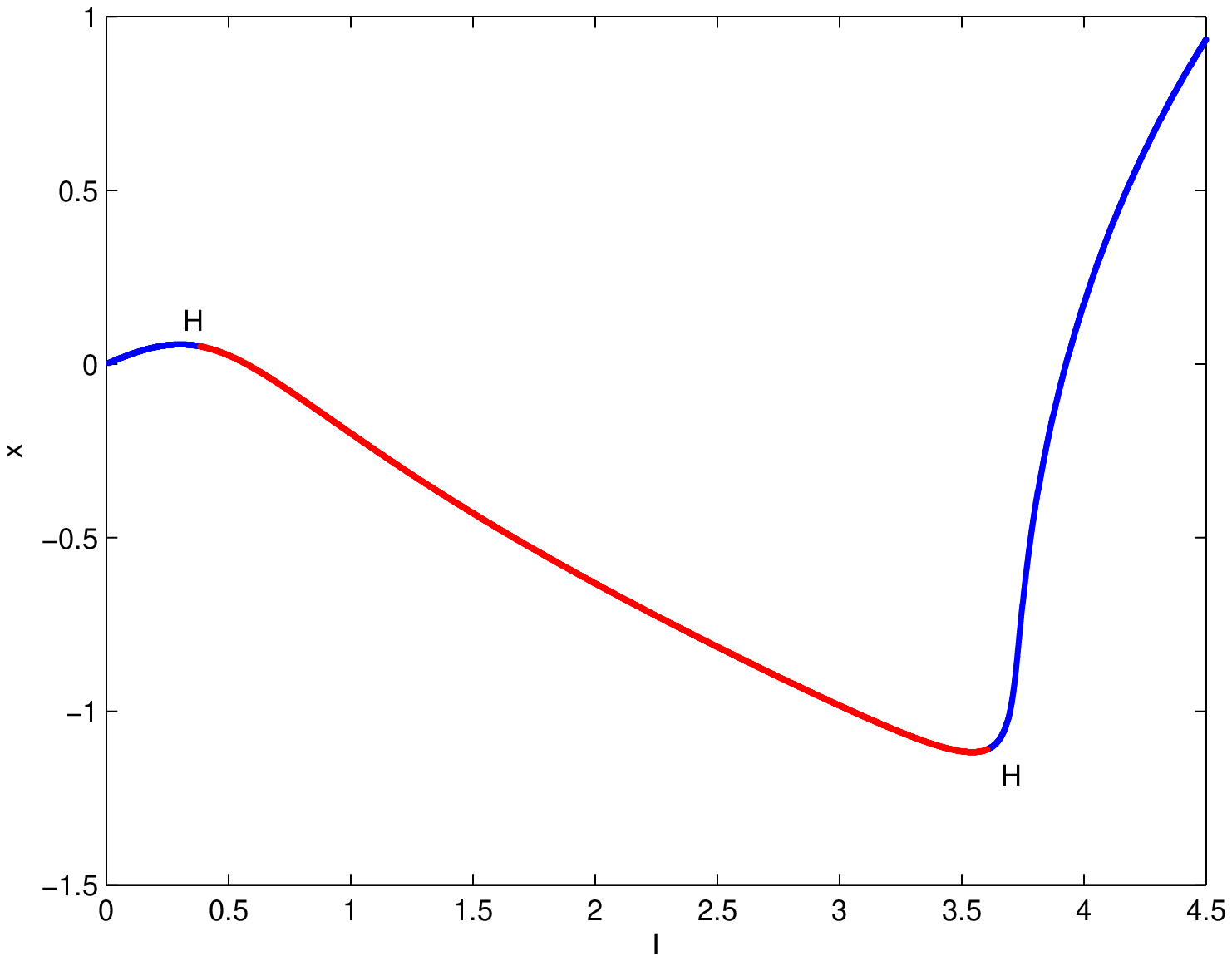}\\
   (e)
 \end{minipage}
\qquad
 \begin{minipage}[r]{5.5cm}
    \centering
    \includegraphics[bb=81 242 507 592,width=55mm,clip]{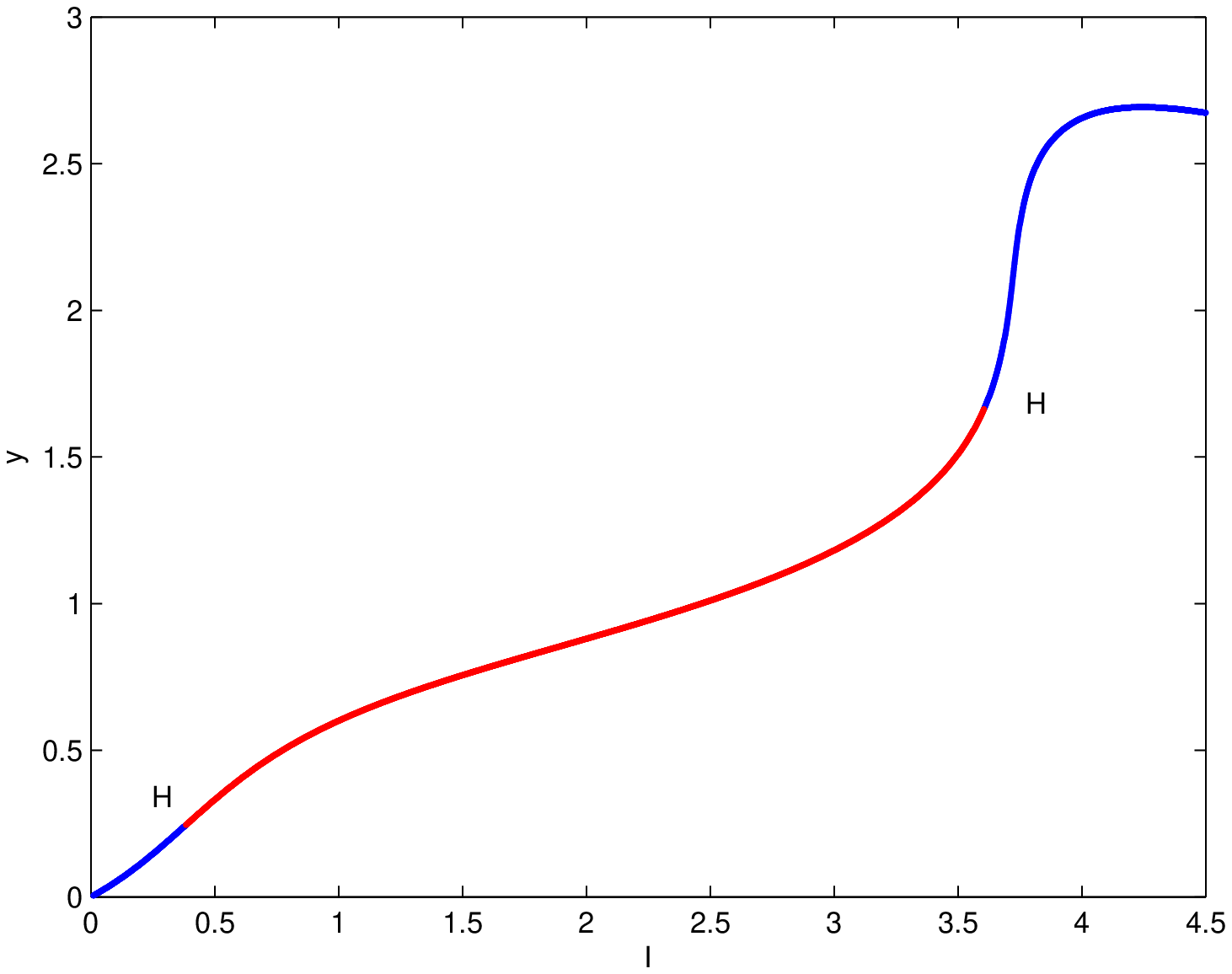}\\
    (f)
 \end{minipage}
 \caption{Stable (in blue) and unstable (in red) equilibrium
configurations of system \eqref{eq:systcart_input}
as function of the parameter I, for different values of $\Omega$ (in (a) and
(b) $\Omega=0.25$, in (c) and (d) $\Omega=0.75$, in (e) and (f) $\Omega=1.5$).
Here $\alpha$ is set equal to 3. In the left and right columns the x and
y-components are represented,
respectively. The bifurcations are indicated as follows:
 SN= saddle-node, SNLC= saddle-node on limit cycle, H= Hopf.} \label{fig:eqOm}
\end{figure}

By choosing $\alpha=3$, in Figure \ref{fig:eqOm} the equilibrium configurations
and their stability properties as function of $I$ for different fixed values of
$\Omega$ are
represented. Depending on the values of the two parameters
$\Omega$ and $I$, a different number of solutions and therefore
different dynamical behaviors are admissible.

\section{Bifurcation analysis}


In this section we study the bifurcations occurring in system
\eqref{eq:systcart_input} as the parameters $(\Omega,I)$ are
varied.

We have seen in the previous section that in absence of external
input ($I=0$) the system presents two nested periodic solutions,
one stable and one unstable, for every value of $\Omega$. For
small values of $I$, these solutions are mantained, but,
increasing the intensity of the input, they disappear through a
sequence of different both local and global bifurcations.

\subsection{Local bifurcations}

\begin{figure}[t!]
    \centering
   \includegraphics[bb=104 271 476 564,width=60mm,clip]{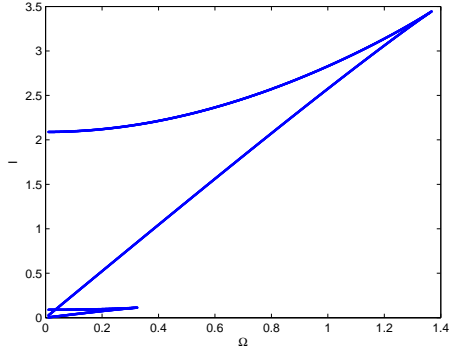}\\
\caption{The locus of parameters where the discriminant of polynomial
\eqref{eq:polyeq} with $\alpha=3$ annihilates. This points out the possible
occurrence of
bifurcations of equilibria
in the system.} \label{fig:discr}
\end{figure}

In order to investigate local bifurcations of equilibria in system
\eqref{eq:systcart_input} we shall look at the polynomial
discriminant of \eqref{eq:polyeq}, and at the linearization of
\eqref{eq:systcart_input} in the neighborhood of the equilibrium
points. A polynomial discriminant is defined as the product of the
squares of the differences of the polynomial roots $s_i$. For a
polynomial of degree $n$ in the form
$$q(x)=a_nx^n+a_{n-1}x^{n-1}+\cdots + a_1x+a_0=0$$
the discriminant is defined as \citep{cohen}
\begin{equation}
D_n=a_n^{2n-2}\prod_{\begin{array}{c} \scriptstyle{i,j}\\ \scriptstyle{ i<j}
\end{array}}^n (s_i-s_j)^2.
\end{equation}
Since the discriminant vanishes in presence of a multiple root,
the values of $\Omega$ and $I$ for which the discriminant $D$ is
equal to zero identify loci of coalescences or births of
solutions, and therefore possible bifurcations. In Figure
\ref{fig:discr} the curves $D=0$ for \eqref{eq:polyeq} in the
$\Omega-I$ plane for $\alpha=3$ are represented.

Now we turn our attention at the linearization of
\eqref{eq:systcart_input}. The Jacobian matrix of system
\eqref{eq:systcomplex_input_risc2}, evaluated in the equilibrium
points $(\overline x, \overline y)$, has the following expression:
\begin{equation}\label{eq:jacob}
 J=\begin{pmatrix}
    g(\overline x,\overline y)+\overline x^2h(\overline x,\overline y) &
\overline x \overline yh(\overline x,\overline y)-\Omega\\
    \overline x\overline yh(\overline x,\overline y)+\Omega &
g(\overline x,\overline y)+\overline y^2h(\overline x,\overline y)
   \end{pmatrix},
\end{equation}
where
\begin{equation}\label{eq:h}
 h(\overline x,\overline y)=\frac{\alpha}{\sqrt{\overline x^2+\overline y^2}}-2.
 \end{equation}
The sign of the real part of the eigenvalues can be determined by
looking at the trace and the determinant of the Jacobian matrix.
In the present case they are given by
\begin{align}
 \tr J & = 2g(\overline x,\overline y)+(\overline x^2+\overline y^2)h(\overline
           x, \overline y)\nonumber\\
        & = 2(-1+\alpha\sqrt{\overline x^2+\overline y^2}-(x^2+y^2))
             +(\overline x^2+\overline y^2)\left(\frac{
             \alpha}{\sqrt{\overline x^2+\overline y^2}}-2\right)\nonumber\\
       & = 2(-1+\alpha\overline \zeta-\overline\zeta^2)+
           \alpha \overline\zeta -2 \overline\zeta^2\nonumber\\
       & = -4\overline \zeta^2 +3\alpha\overline \zeta-2, \label{eq:trace}\\
       & \nonumber\\
 \det J & = (g(\overline x,\overline y)+\overline x^2h(\overline x,\overline y))
            (g(\overline x,\overline y)+\overline y^2h(\overline
            x,\overline y))-(\overline x\overline yh(\overline x,\overline
            y)-\Omega)
            (\overline x\overline yh(\overline x,\overline y)+\Omega)\nonumber\\
        & = (g(\overline x,\overline y))^2+(\overline x^2+\overline
            y^2)g(\overline x,\overline y)h(\overline x,\overline
            y)-\Omega^2\nonumber\\
         & = \frac{\Omega I}{\overline y}-I\overline x h(\overline x,\overline
                y)\nonumber\\
        & = \frac{\Omega I}{\overline y}+\frac{I}{\Omega}\overline y
                \left[-1+\alpha\sqrt{\frac{I}{\Omega} \overline y}
                -\frac{I}{\Omega}\overline y\right]
                \left(\frac{\alpha}{\sqrt{\frac{I}{\Omega} \overline
                y}}-2\right)\nonumber\\
        & =  \frac{I^2}{\overline \zeta^2}+\overline \zeta
                \left[-1+\alpha \overline \zeta - \overline \zeta^2\right]
                \left(\alpha-2\overline \zeta\right), \label{eq:det}
\end{align}
having exploited relations \eqref{eq:g}, \eqref{eq:intersez}, and \eqref{eq:h},
and having introduced the ancillar variable $\overline \zeta=
\sqrt{\frac{I}{\Omega}\overline y}$.

Time-domain simulations show that for certain values of $\Omega$
and $I$ system \eqref{eq:systcart_input} exhibits Hopf
bifurcations. The conditions for the occurrence of this
bifurcation are the following \citep{kuznetsov}:
\begin{equation}\label{eq:condhopf}
 \begin{cases}
   \tr J  =0\\
 \det J  >0.
 \end{cases}
\end{equation}
From \eqref{eq:trace}, condition $\tr J=0$ leads to
\begin{equation}\label{eq:trzero}
-4\overline \zeta^2 +3\alpha\overline \zeta-2=0,
\end{equation}
that is
\begin{equation}\label{eq:zetai}
 \overline \zeta_{1,2}=\frac{3\alpha \pm \sqrt{9\alpha^2-32}}{8}.
\end{equation}
These two solutions exist only for $9\alpha^2
-32\geq0$, that is only for $\alpha\geq \frac{4\sqrt 2}{3}\cong 1.886$. However,
this condition is alway satisfied, since $\alpha>2$ due to
\eqref{eq:cond2cycles}.

From the definition of $\overline \zeta$ and from
\eqref{eq:intersez}, we can find the coordinates of the
equilibrium points $P_1=\left(\overline x_1,\overline y_1\right)$
and $P_2=\left(\overline x_2,\overline y_2\right)$, at which the
trace of the Jacobian annihilates:
\begin{align*}
\overline y_i & =\frac{\Omega}{I} \overline \zeta_i^2\\
                   \overline
x_i & =-\frac{1}{\Omega}\left[-1+\alpha\overline \zeta_i-\overline
\zeta_i^2\right] \qquad i=1,2.
\end{align*}
Substituting \eqref{eq:zetai} into \eqref{eq:polyeq} we obtain the
Hopf bifurcation curves in the parameters plane
\begin{equation}\label{eq:IOmega_hopf}
I^2=I_i^2=\overline \zeta_i^6-2\alpha \overline
\zeta_i^5+(\alpha^2+2)\overline \zeta_i^4-2\alpha\overline
\zeta_i^3+(1+\Omega^2)\overline \zeta_i^2 \qquad i=1,2.
\end{equation}
For instance, with $\alpha=3$ (the case considered in Figure
\ref{fig:eqOm}) we have
\begin{align*}
 \overline \zeta_1=\frac{1}{4} \quad \Rightarrow \quad &
P_1=\left(-\frac{5}{16\Omega},\frac{\Omega}{16I_1}\right)\\
                   & I_1^2=\frac{25+256\Omega^2}{4096}\\
\overline \zeta_2=2 \quad \Rightarrow \quad &
P_2=\left(-\frac{1}{\Omega},\frac{4\Omega}{I_2}\right)\\
                   & I_2^2=4(1+\Omega^2).
\end{align*}

Actually, we still have to discriminate between Hopf bifurcations
and neutral saddles \citep{kuznetsov}, exploiting the condition on
the determinant of the Jacobian matrix. Introducing again the
variable $\overline \zeta_i=\sqrt{\frac{I_i}{\Omega}\overline
y_i}$, and recalling the expression for $I_i$ in
\eqref{eq:IOmega_hopf}, we get
\begin{align*}\label{eq:detzeta}
 \det J_{|_{P_i}} & = \frac{I_i^2}{\overline \zeta_i^2}+\overline \zeta_i
\left[-1+\alpha \overline \zeta_i - \overline \zeta_i^2\right]
\left(\alpha-2\overline \zeta_i\right)\\
&  = \frac{\overline \zeta_i^6-2\alpha
\overline \zeta_i^5+(\alpha^2+2)\overline \zeta_i^4-2\alpha\overline
\zeta_i^3+(1+\Omega^2)\overline \zeta_i^2}{\overline \zeta_i^2}+\zeta_i
\left[-1+\alpha \overline \zeta_i - \overline \zeta_i^2\right]
\left(\alpha-2\overline \zeta_i\right)\\
& = 3\overline \zeta_i^4-5\alpha \overline \zeta_i^3+2(\alpha^2+2)\overline
\zeta_i^2-3\alpha\overline \zeta_i+(1+\Omega^2).
\end{align*}

Requiring the positivity of $\det J$, it yields
\begin{align}
 \det J_{|_{P_i}}>0 \quad \Rightarrow \quad
\Omega^2 & >-3\overline \zeta_i^4+5\alpha \overline
\zeta_i^3-2(\alpha^2+2)\overline \zeta_i^2+3\alpha\overline \zeta_i-1\\
         & = \left(\frac{3}{64}\alpha^3-\frac{1}{4}\alpha\right)\overline
\zeta_i+\frac{1}{4}-\frac{1}{32}\alpha^2,
\nonumber
\end{align}
where we used \eqref{eq:trzero} to reduce the degree of the right
hand side. In particular, substituting the values of $\overline
\zeta_i$ found in \eqref{eq:zetai}, we obtain
\begin{align}
 \det J_{|_{P_1}}>0 \quad \Rightarrow \quad &
\Omega^2>\frac{9}{512}\alpha^4+\frac{1}{4}-\frac{1}{8}\alpha^2-\left(\frac{3}{
512}\alpha^3-\frac{1}{32}\alpha\right)\sqrt{9\alpha^2-32}\\
\det J_{|_{P_2}}>0 \quad \Rightarrow \quad &
\Omega^2>\frac{9}{512}\alpha^4+\frac{1}{4}-\frac{1}{8}\alpha^2+\left(\frac{3}{
512}\alpha^3-\frac{1}{32}\alpha\right)\sqrt{9\alpha^2-32}.\nonumber
\end{align}

\begin{figure}[t!]
    \centering
   \includegraphics[bb=81 242 507 588,width=70mm,clip]{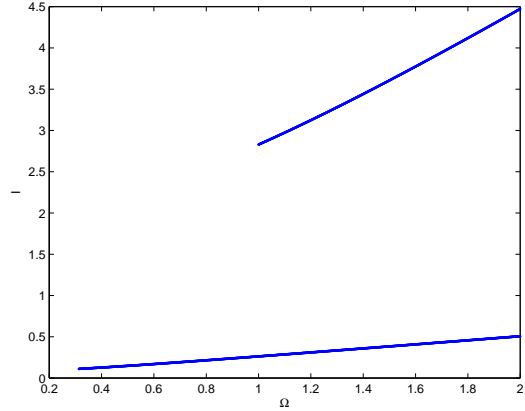}\\
\caption{Curves of Hopf bifurcations in the $\Omega - I$ plane for $\alpha=3$.}
\label{fig:hopf}
\end{figure}

Finally, we can conclude that for $\alpha=3$ in the plane $\Omega - I$ the Hopf
bifurcations
occur if the following conditions are satisfied (see Figure
\ref{fig:hopf}):

\begin{equation}\label{eq:hopf}
 \begin{cases}
  I^2=\frac{25+256\Omega^2}{4096}\\
  \Omega>\frac{5}{16}
 \end{cases}
 \qquad
 \begin{cases}
  I^2=4(1+\Omega^2)\\
  \Omega>1.
 \end{cases}
 \qquad
\end{equation}

It is worth observing that, for the examples shown in Figure \ref{fig:eqOm}, we
have one Hopf bifurcation for $\Omega=0.75$ (in this case only $P_1$ exists) and
two Hopf bifurcations for $\Omega=1.5$. These situations are in perfect
agreement with conditions \eqref{eq:hopf}.

\begin{remark}
The points $P_1$ and $P_2$ found above have a direct relation with
the curves $\Gamma$ and $S$ in \eqref{eq:intersez}. If we consider
$S$ as $S: x=f(y)$, then it is possible to show that $P_1$ and
$P_2$ are the relative maximum and minimum of $S$, respectively.
Furthermore, the values $I_1$ and $I_2$ are precisely the ones for
which these relative maximum and minimum points lie on the
circumference $\Gamma$.
\end{remark}

To determine the stability of the emerging limit cycle, we use
normal form theory. We recall \citep{guckenheimer} that for a
planar system in the form
\begin{equation}
 \begin{pmatrix}
  \dot x\\
  \dot y\\
 \end{pmatrix}=
 \begin{pmatrix}
  0 & -\omega\\
  \omega & 0\\
 \end{pmatrix}
 \begin{pmatrix}
   x\\
   y\\
 \end{pmatrix}+
 \begin{pmatrix}
  F(x,y)\\
  G(x,y)\\
 \end{pmatrix}
\end{equation}
with $F(0)=G(0)=0$ and $DF(0)=DG(0)=0$, we have to evaluate the
following quantity:
\begin{align}\label{eq:a}
 a=\frac{1}{16} & \left(F_{xxx}+F_{xyy}+G_{xxy}+G_{yyy}\right)\\
 &+\frac{1}{16\omega}
\left(F_{xy}(F_{xx}+F_{yy})-G_{xy}(G_{xx}+G_{yy})-F_{xx}G_{xx}+F_{yy}G_{yy}
\right),\nonumber
\end{align}
where all the partial derivatives are computed in $(0,0)$. Thus,
if $a<0$ we can conclude that the Hopf bifurcation is
supercritical, while if $a>0$ is subcritical \citep{guckenheimer}.

In order to proceed with the computation we have to move the
bifurcation points from $P_1$ and $P_2$ to the origin. Let us
introduce the following change of variables:
\begin{equation}
 \begin{cases}
  u=x-\overline x_i\\
  v=y-\overline y_i,\\
 \end{cases}
\end{equation}
where $i=1,2$, depending on the point we are interested in. Thus, our system
\eqref{eq:systcart_input} can be recast as
\begin{equation}
 \begin{pmatrix}
  \dot u\\
  \dot v\\
 \end{pmatrix}=
 \begin{pmatrix}
  0 & -\Omega\\
  \Omega & 0\\
 \end{pmatrix}
 \begin{pmatrix}
   u\\
   v\\
 \end{pmatrix}+
 \begin{pmatrix}
  F(u,v)\\
  G(u,v)\\
 \end{pmatrix},
\end{equation}
where
\begin{align*}
 F(u,v) & =\left(-1+\alpha\sqrt{(u+\overline x_i)^2+(v+\overline y_i)^2}
               -((u+\overline x_i)^2+(v+\overline y_i)^2)\right)(u+\overline
                x_i)-\Omega \overline y_i+I\\
G(u,v) &= \left(-1+\alpha\sqrt{(u+\overline x_i)^2+(v+\overline y_i)^2}
               -((u+\overline x_i)^2+(v+\overline
                y_i)^2)\right)(v+\overline y_i)+\Omega \overline x_i.
 \end{align*}
It is easy to check that the functions $F(u,v)$ and $G(u,v)$
satisfy the conditions above. Furthermore, evaluating their
partial derivatives, it is possible to observe that $G_u=F_v$,
which implies $G_{uu}=F_{uv}$, $G_{uv}=F_{vv}$, and
$G_{uuv}=F_{uvv}$. Expression \eqref{eq:a} reduces to
\begin{equation}\label{eq:a2}
 a=\frac{1}{16} \left(F_{uuu}+2F_{uvv}+G_{vvv}\right).
\end{equation}
Computing these partial derivatives and evaluating them in $(0,0)$, we get
\begin{equation}\label{eq:a3}
 a=\frac{1}{16} \left(-16+\frac{3\alpha}{\sqrt{\overline
x_i^2+\overline y_i^2}}\right).
\end{equation}
Recalling that from \eqref{eq:circum} we have $\overline
x_i^2+\overline y_i^2=\frac{I}{\Omega}\overline y_i=\overline \zeta_i^2$, we
obtain
\begin{equation}\label{eq:a4}
 a=\frac{1}{16} \left(-16+\frac{3\alpha}{\overline \zeta_i}\right)
\end{equation}
and therefore
\begin{equation}\label{eq:aP1}
 a_{|_{P_1}} =\frac{1}{16}
\left(-16+\frac{3\alpha}{\zeta_1}\right) = -1+\frac{3\alpha
}{2(3\alpha-\sqrt{9\alpha^2-32})}>0.
\end{equation}
Therefore, we conclude that for every $\alpha>2$
and $\Omega>\frac{5}{16}$, in $P_1$  we have a subcritical Hopf
bifurcation. Analogously, evaluating the quantity $a$ in $P_2$ we
obtain:
\begin{equation}\label{eq:aP2}
  a_{|_{P_2}} = \frac{1}{16}
\left(-16+\frac{3\alpha}{\zeta_2}\right)  = -1+\frac{3\alpha
}{2(3\alpha+\sqrt{9\alpha^2-32})}<0.
\end{equation}
In this case we expect a supercritical Hopf bifurcation for every
$\alpha>2$ and $\Omega>1$. Both the results have been confirmed by
time-domain numerical simulations (see also Figure \ref{fig:eqOm}).\\

We now turn our attention to saddle-node bifurcations of
equilibria, that are characterized by $\det J=0$, since they
involve the presence of a null eigenvalue \citep{kuznetsov}. Using
\eqref{eq:det} we obtain
\begin{equation}\label{eq:Ifold}
 I^2=-2\overline
\zeta^6+3\alpha \overline \zeta^5-(2+\alpha^2)
\overline \zeta^4+\alpha \overline \zeta^3.
\end{equation}

Because the Jacobian is evaluated in the equilibrium
configurations, we recall that $\overline \zeta$ will be a root of
\eqref{eq:polyeq}.
Thus, substituting \eqref{eq:Ifold} in \eqref{eq:polyeq} we obtain
\begin{equation}\label{eq:condfold}
 \overline \zeta^2\left[3\overline
\zeta^4-5\alpha \overline \zeta^3+2(2+\alpha^2) \overline
\zeta^2-3\alpha \overline \zeta + (1+\Omega^2)\right]=0.
\end{equation}
We conclude that, for any fixed value of the parameter $\Omega$,
we have a saddle-node bifurcation at the equilibrium points that
satisfy:
\begin{equation*}
 3\overline
\zeta^4-5\alpha \overline \zeta^3+2(2+\alpha^2)
\overline \zeta^2-3\alpha \overline \zeta +
(1+\Omega^2)=0
\end{equation*}
for the values of $I$ given by \eqref{eq:Ifold}.
\begin{remark}
The conditions for the occurrence of a saddle-node bifurcation
leads to the same set of curves in the $\Omega - I$ plane of
Figure \ref{fig:discr}. This is due to the fact that a zero
eigenvalue for the Jacobian matrix of a generic dynamical system
$\dot x=f(x)$ implies that the equilibrium point has a
multiplicity equal to two as zero of the function f(x)=0
\citep{kuznetsov}.
\end{remark}

Finally, we consider the case
\begin{equation}\label{eq:condbt}
 \begin{cases}
   \tr J  =0\\
 \det J  =0
 \end{cases}
\end{equation}
that corresponds to a codimension 2 bifurcation, the so-called
Bogdanov-Takens one. In our case, we have previously seen that the
trace of the Jacobian matrix is equal to zero at the points $P_1$
and $P_2$, with $I_i$ given by \eqref{eq:IOmega_hopf}. The further
condition on the determinant of $J$ leads to the following
critical values of $\Omega$:
\begin{align}\label{eq:omega_bt}
 P_1 \quad \Rightarrow \quad &
\Omega_1^2=\frac{9}{512}\alpha^4+\frac{1}{4}-\frac{1}{8}\alpha^2+\left(\frac{3}{
512}\alpha^3-\frac{1}{32}\alpha\right)\sqrt{9\alpha^2-32}\\
P_2 \quad \Rightarrow \quad &
\Omega^2_2=\frac{9}{512}\alpha^4+\frac{1}{4}-\frac{1}{8}\alpha^2+\left(-\frac{3}
{ 512}\alpha^3+\frac{1}{32}\alpha\right)\sqrt{9\alpha^2-32}.\nonumber
\end{align}

In particular for $\alpha=3$ we have
\begin{align*}
 P_1=\left(-\frac{5}{16\Omega},\frac{\Omega}{16I_1}\right) \quad & \text{with}
\quad
I_1=\frac{\sqrt{25+256\Omega^2}}{8}\\
P_2=\left(-\frac{1}{\Omega},\frac{4\Omega}{I_2}\right) \quad & \text{with} \quad
                    I_2=2\sqrt{1+\Omega^2}
\end{align*}
and, exploting \eqref{eq:omega_bt}, we conclude that in
\begin{equation*}
 (\Omega,I)
=\left(\frac{5}{16},\frac{5}{64}\sqrt{2}\right)\qquad\text{and}\qquad
 (\Omega,I)  =\left(1,2\sqrt{2}\right)
\end{equation*}
we have two Bogdanov-Takens bifurcations.

\subsection{Global bifurcations}
Time-domain numerical simulations reveal that, for some couple
$(\Omega,I)$ shown in Figure \ref{fig:discr}, the saddle-node
bifurcations of equilibria are nontrivial, in the sense that they
involve the appearance and disappearance of limit cycles.

\begin{figure}[t!]
    \centering
   \includegraphics[bb=2 162 582 313,width=100mm,clip]{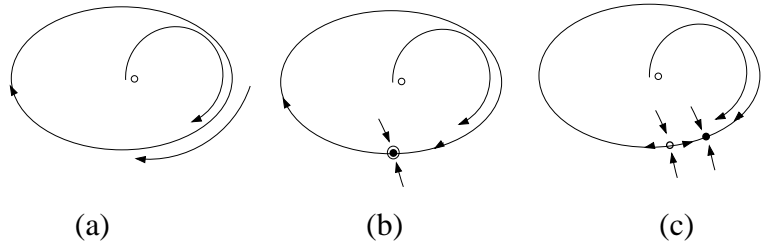}\\
   (a)\hspace{3.3cm} (b)\hspace{3.3cm} (c)\\
\caption{Saddle-node bifurcation on limit cycle (SNLC).}
\label{fig:snlc}
\end{figure}

Firstly introduced and studied in \citep{leontovich}, in system
\eqref{eq:systcart_input} it involves the disappearance of the periodic
solutions. In literature, this bifurcation goes by several names: saddle-node on
limit cycle (SNLC) \citep{hoppen2}, saddle-node on invariant cycle
(SNIC) \citep{izhikevich2}, saddle-node infinite period (SNIPER)
\citep{mccormick} or saddle-node homoclinic bifurcation
\citep{kuznetsov}. In particular, its name saddle-node on
invariant cycle is due to the fact that it is a standard
saddle-node, but occurs on an invariant cycle \citep{izhikevich2}.
Let us suppose that the system exhibits a periodic solution, as in
Figure \ref{fig:snlc} (a). The emergence of a saddle-node (see
Figure \ref{fig:snlc} (b)) coincides with the break of the limit
cycle, that becomes a homoclinic trajectory. As the bifurcation
parameter increases (see Figure \ref{fig:snlc} (c)), the node and
the saddle move away each other and two heteroclinic trajectories
arise to connect the two equilibria.

It is worth remarking that this bifurcation is both global and
local, as it involves a simultaneous collision of equilibria and
manifolds \citep{kuznetsov}. In fact, as we have shown in the
previous section, local analysis describes only the saddle-node
bifurcation of equilibria, missing the disappearance of the limit
cycle.

In general, detecting a homoclinic trajectory is not a simple
task. A possible evidence of a SNLC bifurcation is given by the
period of the limit cycle. In fact, the closer is the parameter to
the critical value, the larger is the period of the corresponding
limit cycle, that tends to infinity approaching the SNLC
bifurcation. Furthermore, the system spends more time near the
place where the saddle-node will appear, in a sort of sense having
a hunch of the future bifurcation.

\begin{figure}[t!]
 \begin{minipage}[l]{6cm}
   \centering
   \includegraphics[bb=82 250 502 588,width=60mm,clip]{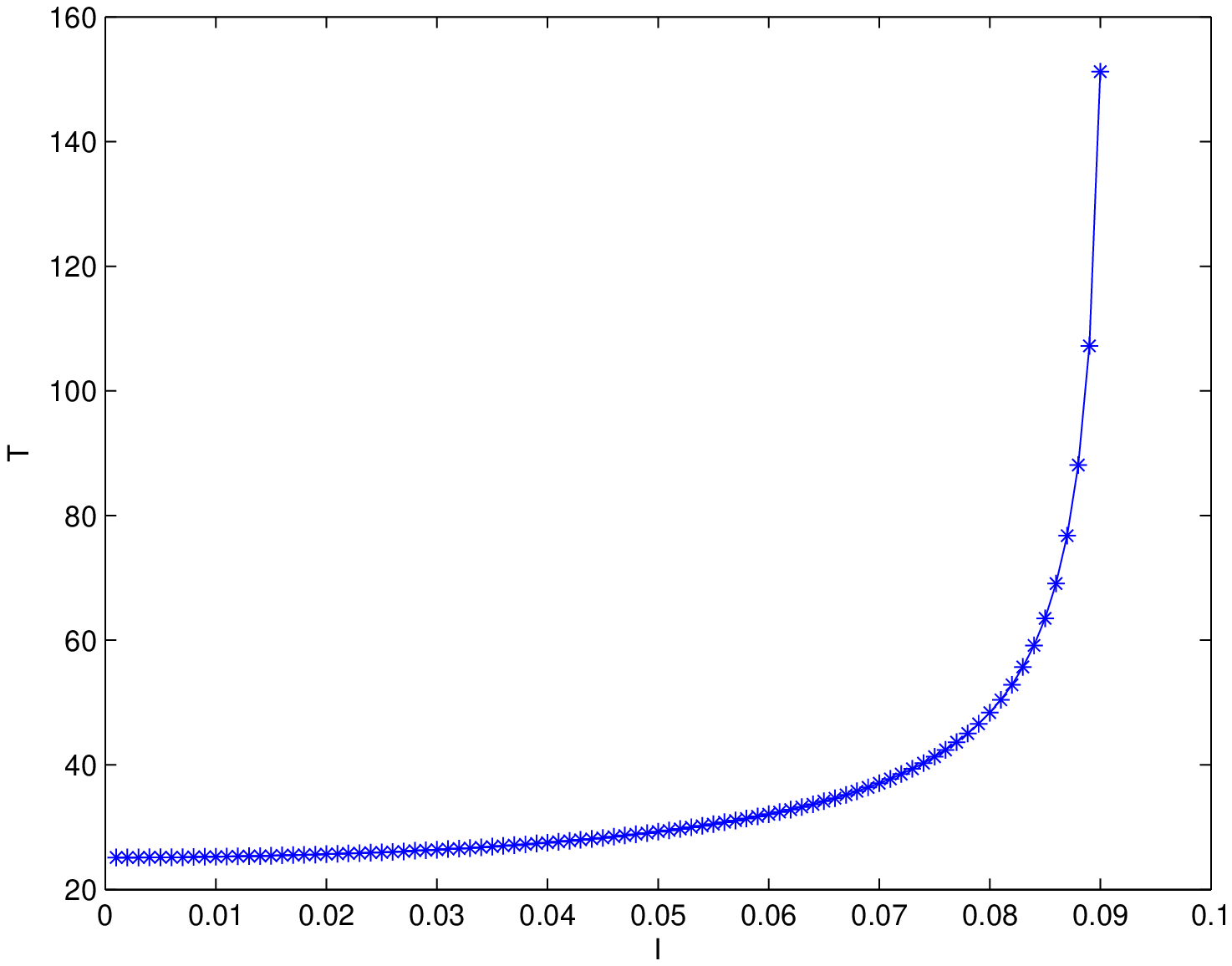}\\
   (a)
 \end{minipage}
\qquad
 \begin{minipage}[r]{6cm}
    \centering
    \includegraphics[bb=82 250 502 588,width=60mm,clip]{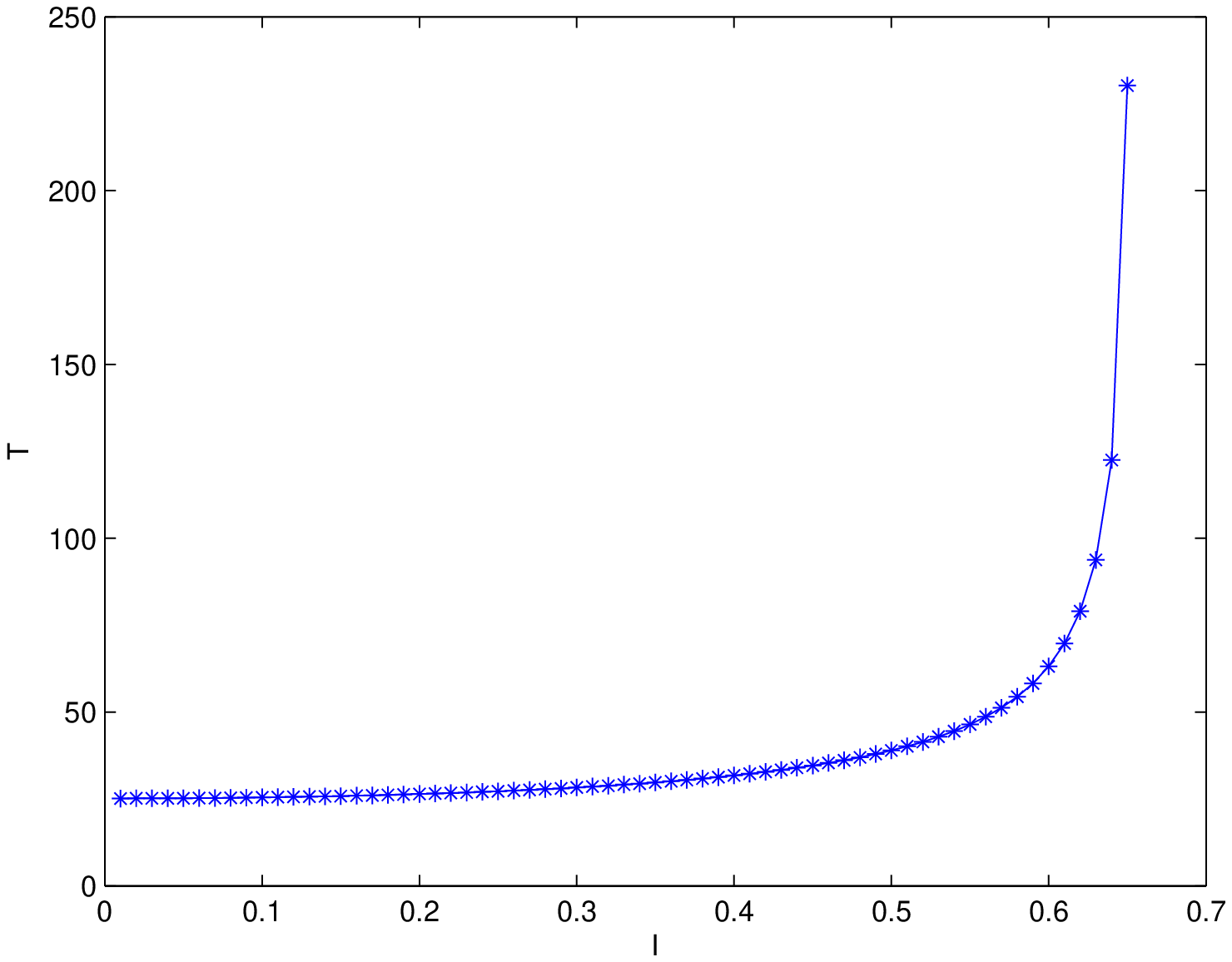}\\
    (b)
 \end{minipage}
\caption{Period $T$ of the unstable (a) and stable (b) limit cycles as function
of the external input $I$, for $\Omega=0.25$. In $I_{c1}=0.091$ and
$I_{c2}=0.653$
the system undergoes two SNLC bifurcations, that entail the
disappearance in succession of the two limit cycles.} \label{fig:T}
\end{figure}

As an example, let us consider the case $\Omega=0.25$, for which
it is possible to see that our system undergoes two SNLC
bifurcations. In Figure \ref{fig:T} the periods of the two (stable
and unstable) limit cycles as function of the external input are
shown. It is worth observing that the two periods tend to
infinity, approaching the bifurcation points $I_{c1}=0.091$ and
$I_{c2}=0.654$ (see also Figures \ref{fig:eqOm}(a-b)),
respectively.

\begin{figure}[t!]
 \begin{minipage}[l]{6cm}
   \centering
   \includegraphics[bb=77 244 515 588,width=60mm,clip]{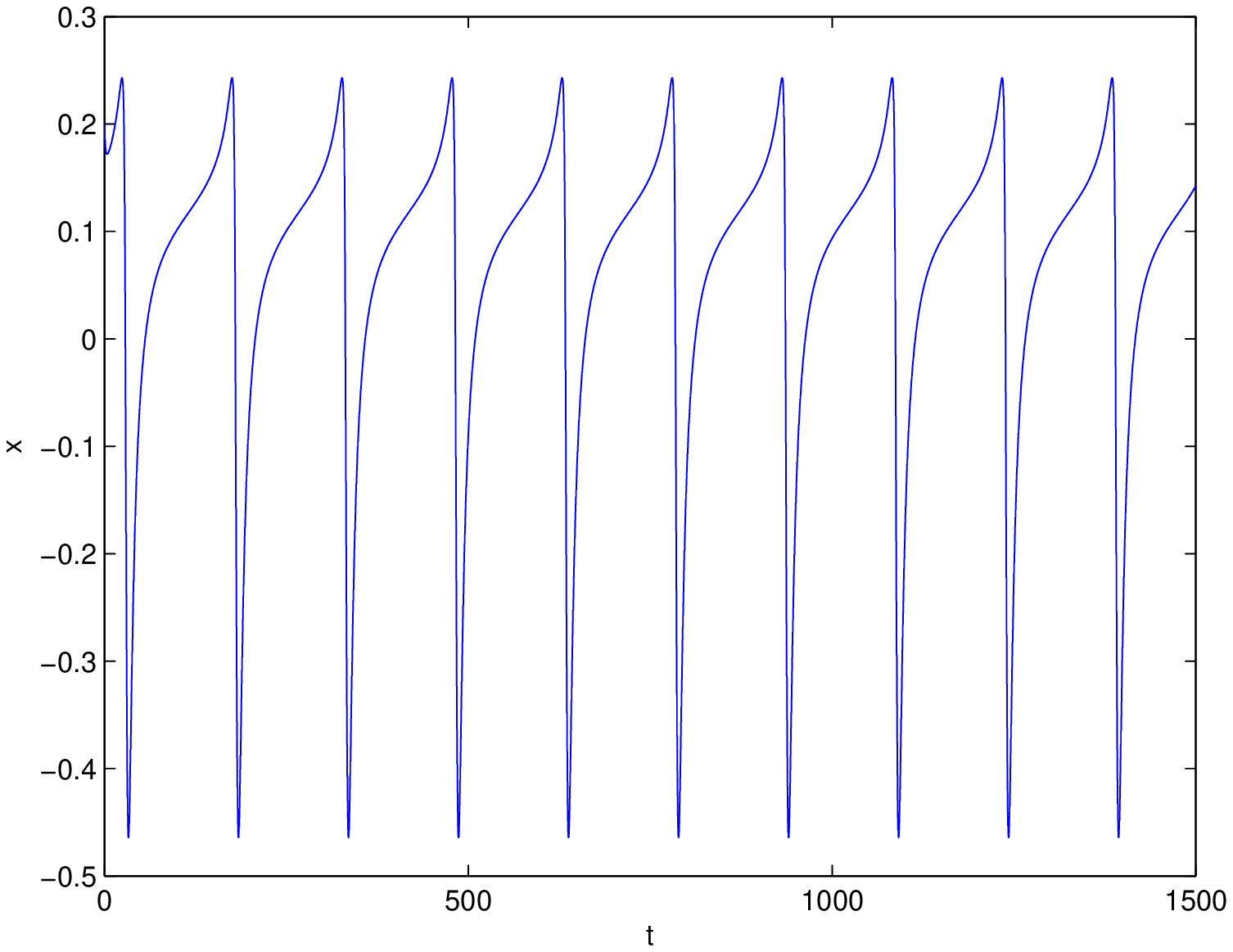}\\
   (a)
 \end{minipage}
\qquad
 \begin{minipage}[r]{6cm}
    \centering
    \includegraphics[bb=77 244 515 585,width=60mm,clip]{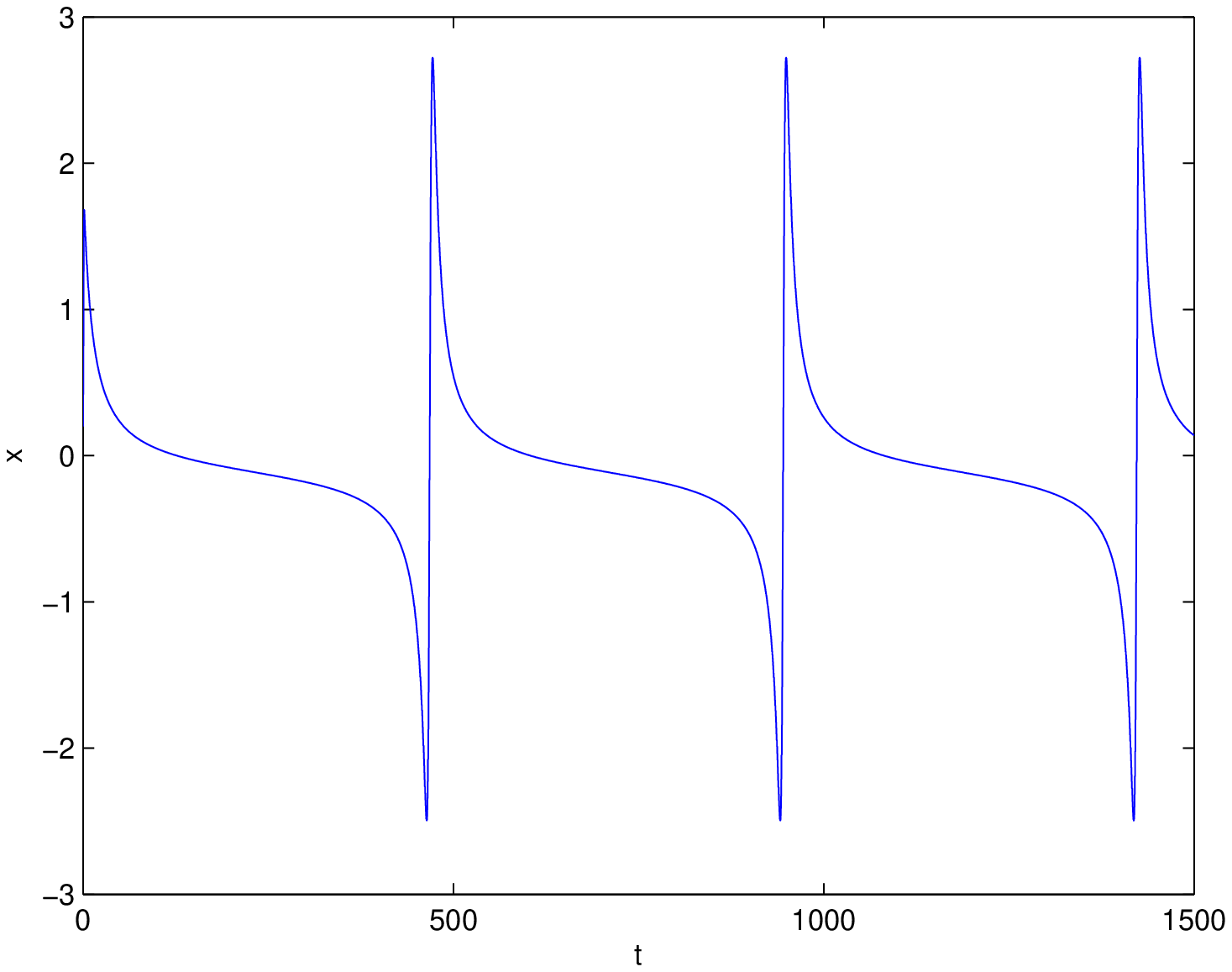}\\
    (b)
 \end{minipage}
\caption{Wave form of the unstable (a) and stable (b) limit cycles
($x$-component), in proximity of SNLC bifurcations ($I=0.090$ and
$I=0.654$, respectively). The frequency $\Omega$ is set equal to $1$.}
\label{fig:x_snlc}
\end{figure}

In Figure \ref{fig:x_snlc} the wave forms for the x-component of
the two limit cycles in proximity of the respective SNLC bifurcations are
represented. Notice the analogy with dynamics of the so-called
relaxation oscillators.

\begin{figure}[t!]
 \begin{minipage}[l]{6cm}
   \centering
   \includegraphics[bb=82 250 502
599,width=60mm,clip]{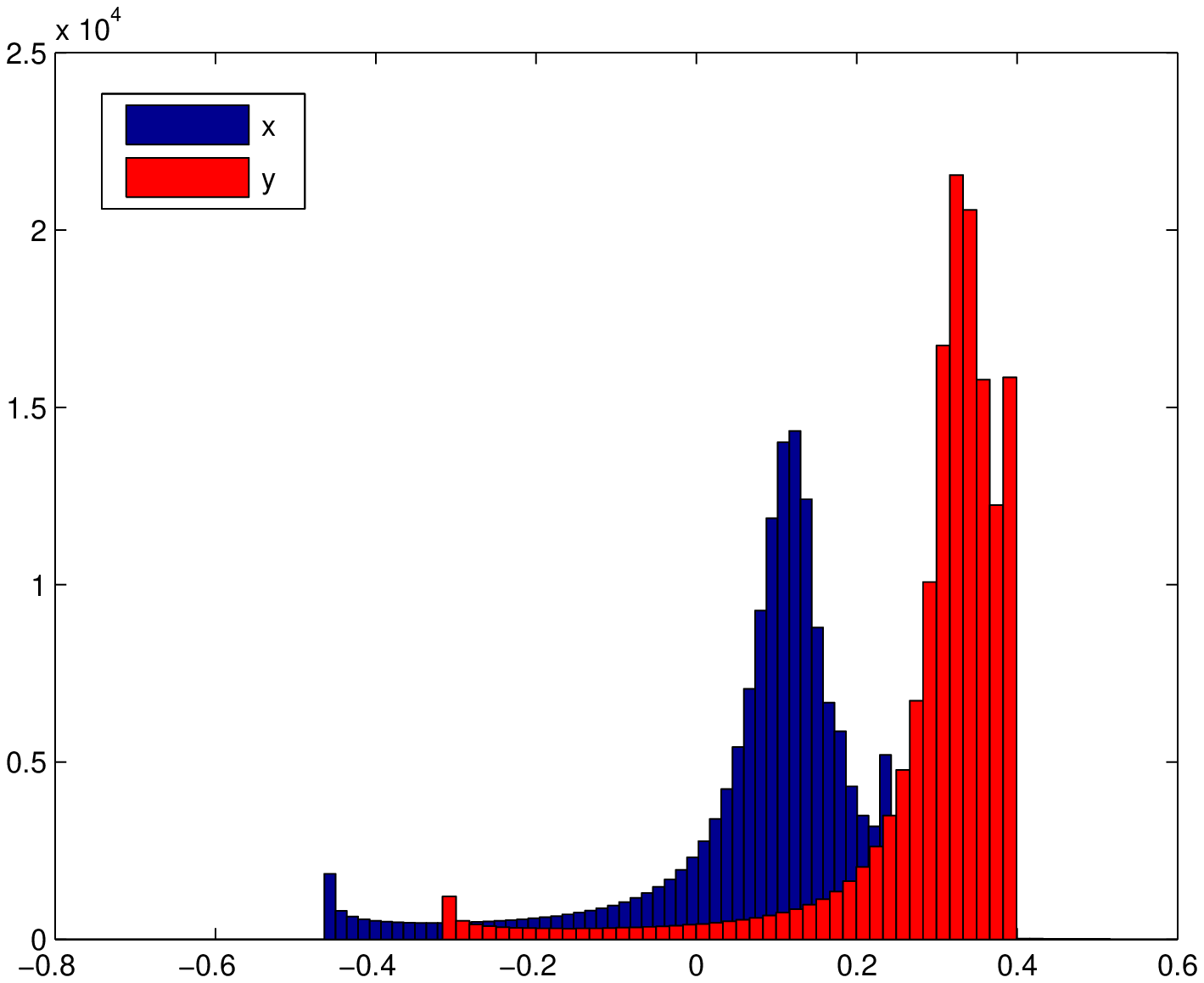}\\
   (a)
 \end{minipage}
\qquad
 \begin{minipage}[r]{6cm}
    \centering
    \includegraphics[bb=82 250 502 599,width=60mm,clip]{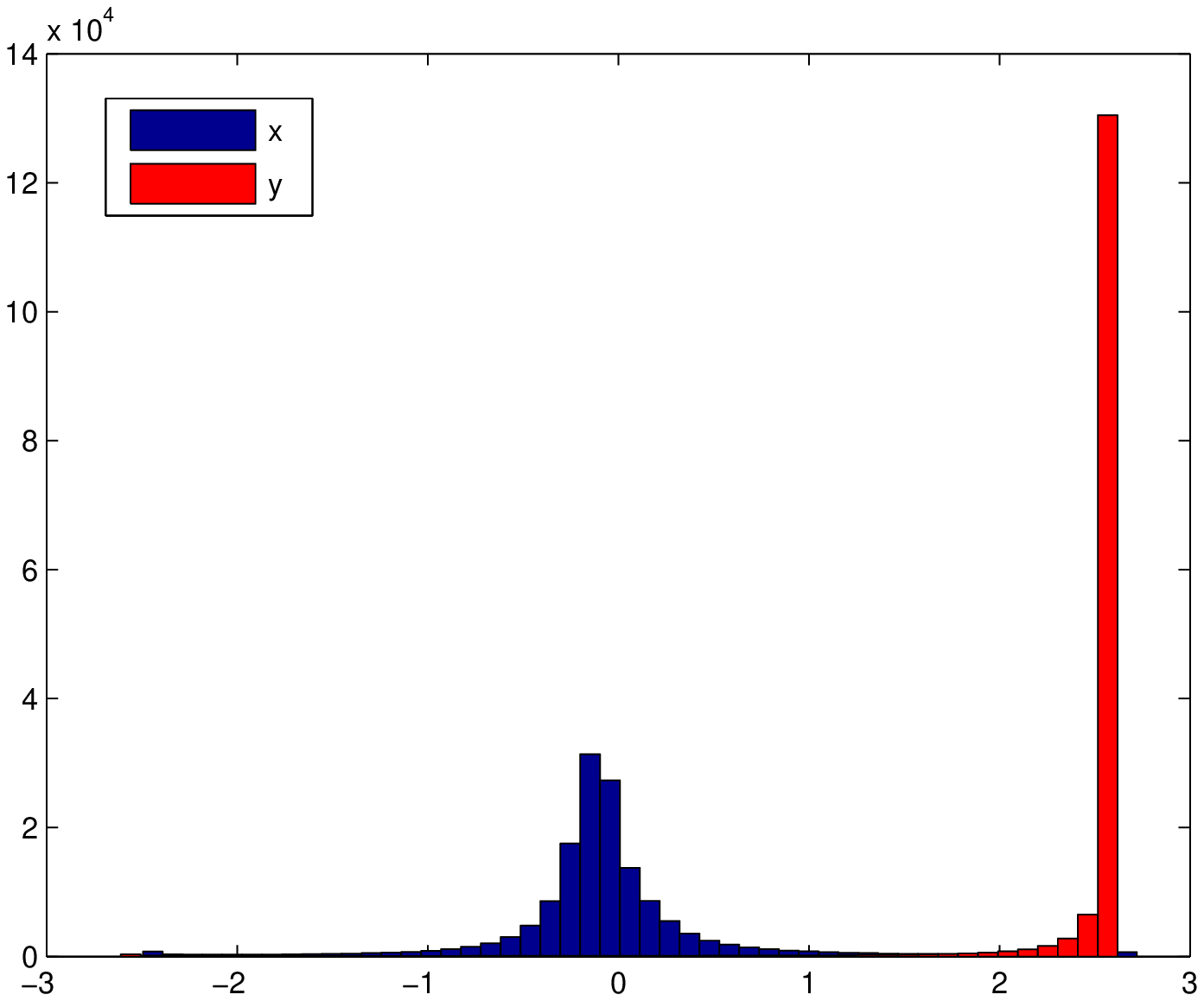}\\
    (b)
 \end{minipage}\\
 \begin{minipage}[l]{6cm}
   \centering
   \includegraphics[bb=82 250 502
599,width=60mm,clip]{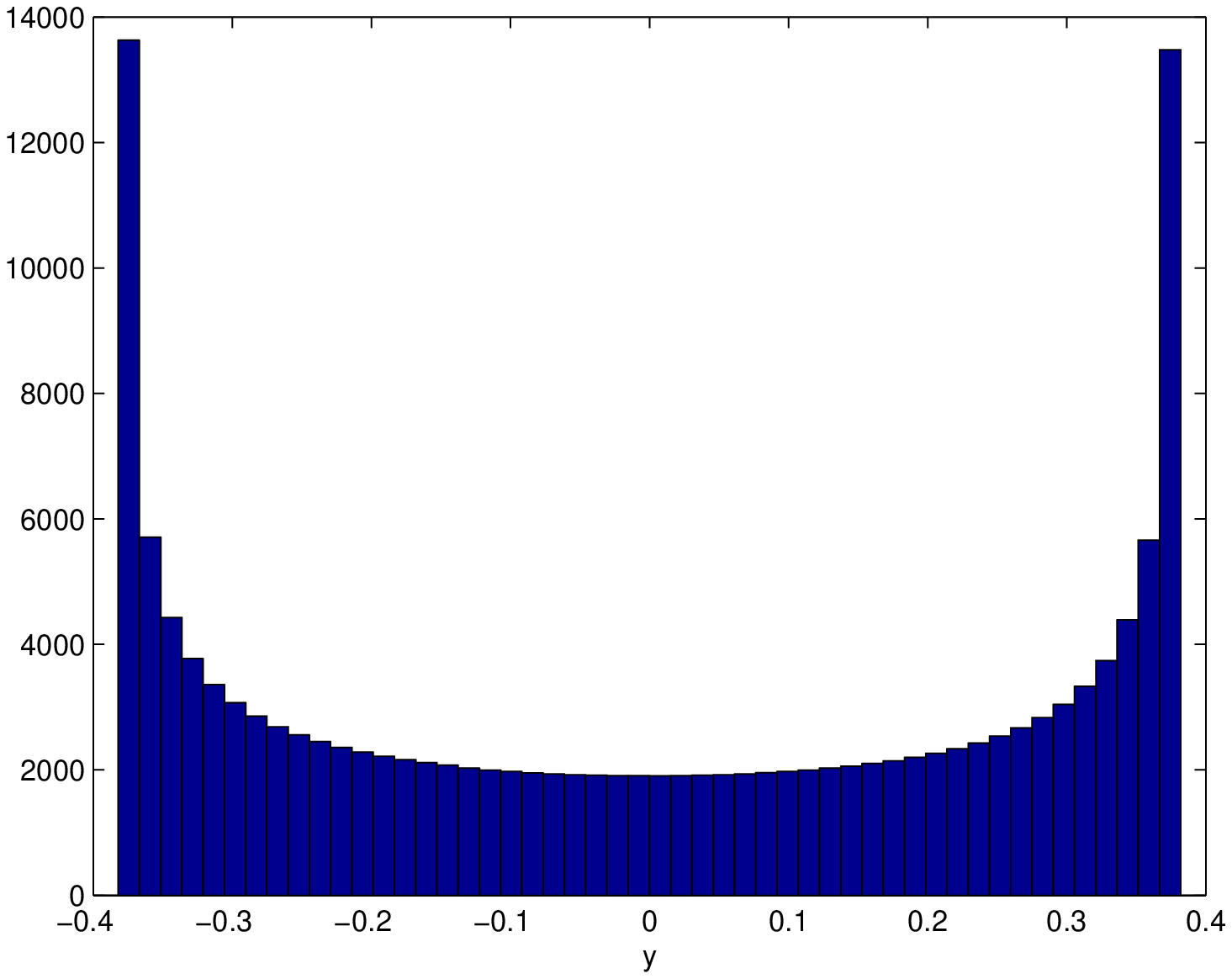}\\
   (c)
 \end{minipage}
\qquad
 \begin{minipage}[r]{6cm}
    \centering
    \includegraphics[bb=82 250 502 599,width=60mm,clip]{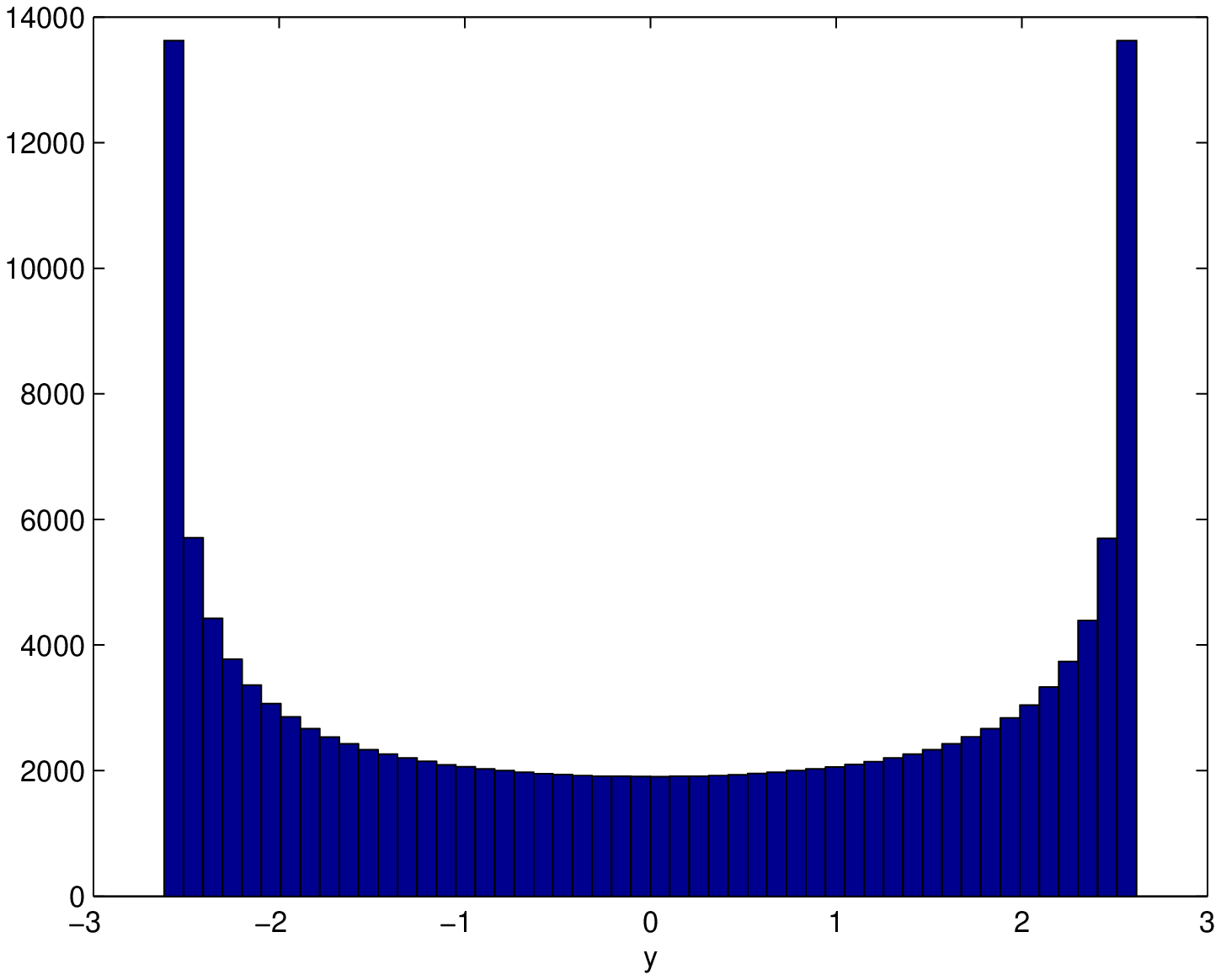}\\
    (d)
 \end{minipage}
\caption{Histogram of the unstable (a-c) and stable (b-d) limit cycle
dynamics near
 the SNLC bifurcation ($I=0.090$ and
$I=0.653$ in (a) and (b), respectively) and far away from the critical values of
the
bifurcation parameter ($I=0$ in both (c) and (d)).}
\label{fig:hist_snlc}
\end{figure}

It is possible to see that the system tends to stay for the
majority of the time around one configuration, that is precisely
where the saddle-node appears at the bifurcation point. In order
to characterize the positions in the phase plane where the system
spends most of the time, a histogram of the dynamics arising from
the simulations in the time-domain ($t_{fin}=1500$ and $\Delta t=0.01$) has been
computed (see Figure \ref{fig:hist_snlc}). Since for
$I_{c1}=0.091$ the saddle-node has coordinates
$SN1=(0.1081,0.3284)$, while for $I_{c2}=0.654$ we have
$SN2=(-0.1368, 2.6066)$, the check with Figure \ref{fig:hist_snlc}
leads to the expected conclusions. The same approach has been
carried out for both the two limit cycles in absence of external
input ($I=0$), to make a comparison with the previous relaxation
oscillator-like behavior. In this case, basically all the states
are uniformly visitated by the sistem. The two peaks at minimum
and maximum values are due to the discretization of the variables
for the histogram computation.

It is worth observing that, in order to obtain Figures \ref{fig:T},
\ref{fig:x_snlc} and \ref{fig:hist_snlc} for the unstable limit cycle, the
numerical simulations have been performed back in time.

\section{Conclusion}

Cyclic Negative Feedback Systems (CNF systems) are one of the most
exploited mathematical frameworks to model phenomena arising in
systems biology, such as cascades of molecular reactions inside
the cell. Since it is well known that these events take place in
different compartments, it seems more appropriate to consider
networks of diffusively coupled CNF systems. In addition, the
effect of an external agent, such as the light, the temperature or
the synthesis of the substrate, is suitably modeled as a constant
external term. Unfortunately, due to the huge number of factors
involved in such mechanisms, networks of coupled CNFs systems can
be high-dimensional, and therefore their dynamics can be difficult
to fully characterize.

Inspired by the complex spatio-temporal patterns found in arrays
of diffusively coupled Cyclic Negative Feedback systems (CNF
systems), we have considered the case of radial isochron clocks
that exhibit the coexistence of different stable attractors, as
well as CNF systems. In fact, these systems present the same
qualitative behavior of CNF systems, but they are quite easy to
handle since they are of lower order. In particular, we have dealt
with systems that present a hard excitation behavior, i.e., that
display at the same time a stable equilibrium point and a stable
limit cycle.

As a first step, we have focused on the effect of a constant
external input on a single radial isochron clock with hard
excitation. We have carried out a characterization of local and
global bifurcations, and in particular, we have detected the
occurrence of saddle-node on limit cycle bifurcations. It is
interesting to notice that, in presence of such bifurcations, the
system exhibits a relaxation oscillator-like dynamics.

Once we have completely analyzed the dynamical behavior of a single system in
presence of a constant external input,
we are now interested in considering networks of
diffusively coupled radial isochron clocks, in order to investigate the complex
dynamics that may arise due to the additional effect of coupling.

\section*{Acknowledgment}
This work was partially supported by the CRT Foundation. L. Ponta acknowledges
the Istituto Superiore Mario Boella for financial support. V. Lanza and M.
Bonnin acknowledge the Istituto Superiore Mario Boella and the regional
government
of Piedmont for financial support.

\bibliographystyle{model5-names}
\bibliography{biblio}







\end{document}